\documentclass[12pt]{amsart}
\usepackage{epsfig}
\usepackage{amsmath}
\usepackage{amsfonts}
\usepackage{amssymb}
\usepackage{epigraph, fancyhdr}

\def\D{\mathcal D}
\def\F{\mathcal F}
\def\K{\mathcal K}

\def\J{\mathcal J}
\def\L{\mathcal L} 
\def\H{\mathcal H}
\def\O{\mathcal O}
\def\T{\mathcal T}

\def\1{\mathbf 1}
\def\M{{\overline{\mathcal M}}}

\def\QQ{\mathbb Q}

\def\CC{\mathbb C}
\def\RR{\mathbb R}

\def\Res{\operatorname{Res}}

\def\hat{\widehat}

\def\p{\partial}
\def\a{\alpha}

\def\f{{\mathbf f}}
\def\g{{\mathbf g}}
\def\t{{\mathbf t}}
\def\x{{\mathbf x}}

\def\la{\lambda}
\def\ge{\epsilon}
\def\gL{\Lambda}

\def\m{{\mathfrak m}}

\def\lan{\langle}
\def\ran{\rangle}

\def\ev{\operatorname{ev}}

\def\td{\operatorname{td}}
\def\ch{\operatorname{ch}}
\def\qch{\operatorname{qch}}

\def\fake{\operatorname{fake}}

\def\z{\zeta}

\renewcommand{\Delta}{\triangle}
\makeatletter
\DeclareFontFamily{OMX}{MnSymbolE}{}
\DeclareSymbolFont{MnLargeSymbols}{OMX}{MnSymbolE}{m}{n}
\SetSymbolFont{MnLargeSymbols}{bold}{OMX}{MnSymbolE}{b}{n}
\DeclareFontShape{OMX}{MnSymbolE}{m}{n}{
    <-6>  MnSymbolE5
   <6-7>  MnSymbolE6
   <7-8>  MnSymbolE7
   <8-9>  MnSymbolE8
   <9-10> MnSymbolE9
  <10-12> MnSymbolE10
  <12->   MnSymbolE12
}{}
\DeclareFontShape{OMX}{MnSymbolE}{b}{n}{
    <-6>  MnSymbolE-Bold5
   <6-7>  MnSymbolE-Bold6
   <7-8>  MnSymbolE-Bold7
   <8-9>  MnSymbolE-Bold8
   <9-10> MnSymbolE-Bold9
  <10-12> MnSymbolE-Bold10
  <12->   MnSymbolE-Bold12
}{}

\let\llangle\@undefined
\let\rrangle\@undefined
\DeclareMathDelimiter{\llan}{\mathopen}%
                     {MnLargeSymbols}{'164}{MnLargeSymbols}{'164}
\DeclareMathDelimiter{\rran}{\mathclose}%
                     {MnLargeSymbols}{'171}{MnLargeSymbols}{'171}
\makeatother

\title[Explicit reconstruction]
      {Permutation-equivariant \\ quantum K-theory VIII. \\
      Explicit reconstruction}

\author[A. Givental]{Alexander GIVENTAL}
\thanks{This material is based upon work supported by the National 
Science Foundation under Grant DMS-1007164, and by the IBS Center for Geometry 
and Physics, POSTECH, Korea.} 

\date{August 8, 2015}

\begin{document}

\begin{abstract}

In Part VII, we proved that the range $\L_X$ of the big J-function in permutation-equivariant genus-0 quantum K-theory is an overruled cone, and gave its adelic characterization. Here we show that the ruling spaces are $D_q$-modules in Novikov's variables, and moreover, that the whole cone $\L_X$ is invariant under a large group of symmetries defined in terms of $q$-difference operators. We employ this for the explicit reconstruction of $\L_X$ from one point on it, and apply the result to toric $X$, when such a point is given by the $q$-hypergeometric function.         

\end{abstract}

\maketitle

\section*{Adelic characterization}

We begin where we left in Part VII: at a description of the range $\L \subset \K$ in the space $\K$ of $K^0(X)\otimes \gL$-value rational functions of $q$ of the J-function of permutation-equivariant quantum K-theory of a given K\"ahler target space $X$:
\[ \J := 1-q+\t(q)+\sum_{\a}\phi_{\a}\sum_{n,d} Q^d \lan \frac{\phi^{\a}}{1-qL}; \t(L),\dots,\t(L)\ran_{0,1+n,d}^{S_n}.\]
We proved that $\L$ is an overruled cone, i.e. it is swept by a family
of certain $\gL [q,q^{-1}]$-modules, called {\em ruling spaces}: 
\[ \L = \bigcup_{t\in \gL_{+}} (1-q)\, S(q)^{-1}_t\, \K_{+},\]
where $S_t$ is a certain family of ``matrix'' functions rational in $q$,
whose construction we are not going to remind here. Let us recall the {\em adelic characterization} of $\L$, which will be our main technical tool.

It is given in terms of the overruled cone $\L^{fake}\subset \hat{K}$ in the space of vector-valued Laurent series in $q-1$, describing the range of the J-function in {\em fake} quantum K-theory of $X$:
\[ \L^{fake} = \bigcup_{t\in \gL} (1-q)T_t, \ \ T_t = S^{fake}(q)_t^{-1}\hat{\K}_{+}.\]
Here $S^{fake}_t$ is some ``matrix'' Laurent series, $T_t$ is a tangent space to $\L^{fake}$, containing $(1-q)T_t$, and tangent to $\L^{fake}$ at all points of the ruling space $(1-q) T_t$.

According to the last section of Part VII, {\em a rational function $\f \in \K$ lies in $\L$ if and only if its Laurent series expansions $\f_{(\z)}$ near $q=1/\z$ satisfy the following three conditions:

 (i) $\f_{(1)} \in \L^{fake}$;

 (ii) when $\z\neq 0,1,\infty$ is a primitive $m$th root of unity,
  \[ \f_{(\z)}(q^{1/m}/\z) \in \L^{(\z)}_t ,\]
a certain subspace in $\hat{K}$, determined by the tangent space $T_t$ to $\L^{fake}$ at the point $\f_{(1)}$; 
    
(iii) when $\z \neq 0,\infty$ is not a root of unity, $\f_{(\z)}$ is a power series in $q-1/\z$, i.e. $\f$ has no pole at $q=1/\z$.}

The subspace $\L^{(\z)}_t$ is described as $\nabla_{\z} \Psi^m(T_t)\otimes_{\Psi^m(\gL)}\gL$, where the Adams operation $\Psi^m$ acts by $\Psi^m(q)=q^m$ and naturally on the $\lambda$-algebra $K^0(X)\otimes \gL$, and
$\nabla_{\z}$ is the operator of multiplication by
\[ e^{\textstyle \sum_{k>0}\left(  \frac{\Psi^k(T^*_X)}{k(1-\zeta^{-k}q^{k/m})}-\frac{\Psi^{km}(T_X^*)}{k(1-q^{km})}\right)}. \]

In its turn, the cone $\L^{fake}\subset \hat{K}$ (and hence its tangent spaces
$T_t$) can be expressed in terms of the cone $\L^H\subset \H$, describing the range of cohomological J-function in the space $\H$ of Laurent series in one indeterminate $z$ with coefficients in $H^{even}(X)\otimes \gL$. Namely, according to the Hirzebruch--Riemann--Roch theorem \cite{Co} in fake quantum K-theory,    
\[ \qch (\L^{fake}) = \Delta \, \L^H,\]
where the {\em quantum Chern character} $\qch: \hat{K}\to \H$ acts by
$\qch q = e^z$, and by the natural Chern character $\ch: K^0(X)\otimes \gL \to H^{even}(X)\otimes \gL$ on the vector coefficients, while $\Delta$ acts as the multiplication in the classical cohomology of $X$ by the {\em Euler--Maclaurin asymptotics} (see \cite{CGi, Co, GiTo}) of the infinite product:
\[ \Delta \sim \prod_{r=1}^{\infty} \td (T_X \otimes q^{-r}) .\]
Using all these descriptions, we are going on explore how the {\em string and divisor equations} of quantum cohomology theory manifest in the genus-0 permutation-equivariant quantum K-theory.

\section*{Divisor equations and $\D_q$-modules}

Let $p_1,\dots, p_K$ be a basis in $H^2(X,\RR)$ consisting if integer numerically effective classes, and let $Q^d=Q_1^{d_1}\cdots Q_K^{d_K}$, where $d_i=p_i(d)$, represent degree-$d$ holomorphic curves in the Novikov ring.
We remind that the Novikov variables are included into the ground $\lambda$-algebra $\gL$. 

The loop space $\H$ of Laurent series in $z$ with vector coefficients in $H^{even}(X)\otimes \gL$ is equipped with the structure of a module over the algebra $\D$ of differential operators in the Novikov variables, so that
$Q_i$ acts as multiplication by $Q_i$, and $Q_i\p_{Q_i}$ acts as $zQ_i\p_{Q_i}-p_i$. The divisor equations in quantum cohomology theory imply (see e.g. \cite{GiE}), that

\medskip

{\em linear vector fields $\f \mapsto (Q_i\p_{Q_i}-p_i/z)\f$ in $\H$ are tangent to $\L^H\subset \H$.} 

\medskip

In follows that {\em the ruling spaces (as well as tangent spaces) of $\L^H$ are $\D$-modules, i.e. are invariant with respect to each differential operator
$D(Q, zQ\p_Q-p,z)$, and moreover the flow $\epsilon \mapsto e^{\epsilon D/z}$ of the vector field
  $\f\mapsto D\f/z$ preserves $\L^H$.}

Indeed, for $\f\in \L^H$, the vector $(Q_i\p_{Q_i}-p_i)\f $ lies in $T_{|f}\L^H$, and hence $(zO_i\p_{Q_i}-p_i)\f$ lies in the same ruling space $zT_{\f}\L^H$ as $\f$ does. Therefore so does $D\f$, and hence $D\f/z \in T_{\f}\L^H$, i.e. the vector field $\f \mapsto D\f/z$ is tangent to $\L^H$.   

Note that the operator $\Delta$ relating $\L^H$ and $\L^{\fake}$ involves multiplication in the {\em commutative} classical cohomology algebra $H^{even}(X)$, but does {\em not} involve Novikov's variables. Consequently,
{\em the tangent and ruling spaces of $\qch(\L^{fake})$ are $D$-modules} too,
and moreover, {\em the flows $\epsilon \mapsto e^{\epsilon D/z}$ preserve $\qch(\L^{fake})$.}

We equip the space $\K$ of vector-valued rational functions of $q$ with the structure of a module over the algebra $\D_q$ of finite difference operators.
It is generated (over the algebra of Laurent polynomials in $q$) by multiplication operators, acting as multiplications by $Q_i$, and translation operators, acting as $P_iq^{Q_i\p_{Q_i}}$, where $P_i$ is the multiplication in $K^0(X)$ by the line bundle with the Chern character $\ch P_i =e^{-p_i}$. 

\medskip

{\tt Proposition} (cf. \cite{GiTo, GiE}). {\em The ruling spaces of the overruled cone $\L\subset \K$  of permutation-equivariant quantum K-theory is are $\D_q$-modules.}
  
\medskip

{\tt Proof.} If $\f \in \L$, it passes the tests (i),(ii),(iii) of adelic characterization. We need to show that $\g:=P_iq^{Q_i\p_{Q_i}}\f$, which obviously lies in $\K$, also passes the tests (and with the same $t\in K^0(X)\otimes \gL_{+}$). This is obvious for test (iii), and is true about test (i) because of the above $\D$-module (and hence $\D_q$-module) property of the ruling spaces $(1-q)T_t$ of $\L^{fake}$. To verify test (ii), we write:
\[ \g_{(\z)} (q^{1/m}/\z) = P_i (q^{1/m})^{Q_i\p_{Q_i}}\z^{-Q_i\p_{Q_i}}\f_{(\z)}(q^{1/m}/\z).\]
First, note that the operator $\nabla_{\z}$ relating $\L^{(\z)}$ with $\Psi^m(T_t)$ does not involve Novikov's variables and commutes with $\D_q$. 
Next, let us elucidate the notation $\Psi^m(T_t)\otimes_{\Psi^m(\gL)}\gL$. In fact the space so indicated consists of linear combinations  $\sum_a \la_a(Q,q) \Psi^m(\f_a)$, where $\f_a \in T_t$, and $\la_a \in \gL[[q-1]]$.
We have the following commutation relations:
\begin{align*} P_i\Psi^m = \Psi^m P_i^{1/m}, \ \ & (q^{1/m})^{Q_i\p_{Q_i}}\Psi^m=q^{Q_i^m\p_{Q_i^m}}\Psi^m=\Psi^m(q^{1/m})^{Q_i\p_{Q_i}}, \\ & \z^{Q_i\p_{Q_i}}\Psi^m=\z^{mQ_i^m\p_{Q_i^m}}\Psi^m=\Psi^m.\end{align*}
Therefore
\begin{align*} P_i (q^{1/m})^{Q_i\p_{Q_i}}\z^{-Q_i\p_{Q_i}} & \left(\sum_a \la_a\Psi^m(\f_a)\right) = \\
\sum_a & \left(q^{1/m}/\z)^{Q_i\p_{Q_i}Q_i}\la_a\right)\, \Psi^m\left( P^{1/m}(q^{1/m})^{Q_i\p_{Q_i}} \f_a\right),\end{align*}
 which lies in $\Psi^m(T_t)\otimes_{\Psi^m(\gL)}\gL$ since $T_t$ is invariant under the operator $P^{1/m}(q^{1/m})^{Q_i\p_{Q_i}}=e^{(zQ_i\p_{Q_i}-p_i)/m}$. \qed

 \medskip

 Let $D(Pq^{Q\p_{Q}},q)$ be a constant coefficient finite difference operator, by which we mean a Laurent polynomial expression in translation operators $P_iq^{Q_i\p_{Q_i}}$, and maybe $q$, with coefficients from $\gL$ {\em independent of $Q$}. We assume below that $\epsilon \in \gL_{+}$ to assure $\epsilon$-adic convergence of infinite sums.

 \medskip

 {\tt Theorem 1.} {\em The operator
 \[ e^{\textstyle \sum_{k>0} \Psi^k(\epsilon D(Pq^{kQ\p_{Q}}, q))/k(1-q^k)} \]
 preserves $\L\subset \K$.}     

 \medskip

 {\tt Proof.} We show that if $(1-q) \f$ passes tests (i), (ii), (iii) of
 the adelic characterization of $\L$, then $(1-q)\g$, where 
 \[ \g := e^{\textstyle \sum_{k>0} \Psi^k(\epsilon D(Pq^{kQ\p_{Q}}, q))/k(1-q^k)}\, \f,\]
 also does.

 (i) Suppose $(1-q)\f_{(1)}$ lies in the ruling space $(1-q) T_t \subset \L^{fake}$.
 Note that the exponent $\sum_{k>0} \Psi^k(\epsilon D(Pq^{kQ\p_{Q}}, q))/k(1-q^k)$ has first order pole at $q=1$. According to the discussion above
the flow defined by such an operator on $\hat{K}$ preserves $\L^{fake}$, and therefore maps its tangent spaces to tangent spaces, and ruling spaces to ruling spaces, and moreover, the operators regular at $q=1$ preserve each ruling and tangent space. It follows that $(1-q)\g_{(1)}\in (1-q)T_{t'}\subset \L^{fake}$, where
 \[ T_{t'}:=e^{ \textstyle \sum_{k>0} \Psi^k(\epsilon D(Pq^{kQ\p_{Q}}, 1))/k^2(1-q)}T_t.\]

 (ii) We have
 \[ \Psi^m(T_{t'})=e^{\textstyle \sum_{k>0} \Psi^{mk}(\epsilon D(Pq^{kQ\p_{Q}}, 1))/k^2(1-q^m)}\Psi^m(T_t).\]
 On the other hand, for a primitive $m$th root of unity $\z$, 
 \begin{align*} \g_{(\z)}(q^{1/m}/\z)&= e^{\sum_{k>0} \Psi^k(\epsilon D(P(q^{1/m}/\z)^{kQ\p_{Q}}, q^{1\m}/\z))/k(1-q^{k/m}/\z^k)} \f_{(\z)}(q^{\1/m}/\z) \\
   = e^{\textstyle A} &\, e^{\textstyle \sum_{l>0} \Psi^{ml}(\epsilon D(Pq^{lQ\p_{Q}}, 1))/ml(1-q^l)} \f_{(\z)}(q^{1/m}/\z),\end{align*}
 where $A$ is some operator regular at $q=1$. It comes out of refactoring $e^{A+B/(1-q)}$, where $A$ and $B$ are regular at $q=1$, as $e^A e^{B/(1-q)}$. We use here the fact that the operators $A$ and $B$ have constant coefficients, and hence commute. 

 Note that the exponents $\sum_{k>0} \Psi^{mk}(\epsilon D(Q, Pq^{kQ\p_{Q}}, 1))/k^2(1-q^m)$ and $\sum_{l>0} \Psi^{ml}(\epsilon D(Q, Pq^{lQ\p_{Q}}, 1))/ml(1-q^l)$
 agree modulo terms regular at $q=1$ (which, again, commute with the singular terms).  
 Since we are given that
 \[ \f_{(\z)}(q^{1/m}/\z) \in \nabla_{\z}\, \Psi^m(T_t)\otimes_{\Psi^m(\gL)}\gL,\]
 and since $\nabla_{\z}$ commutes with $\D_q$, we conclude (using the refactoring again), that 
\[ \g_{\z}(q^{1/m}/\z)\in \nabla_{\z} \, \Psi^m(T_{t'})\otimes_{\Psi^m(\gL)}\gL.\]
Note that the exponent in $e^A$ involves translations $P_iq^{Q_i\p_{Q_i}}$ as well as $\z^{-Q_i\p_{Q_i}}$, and so it is important, that (as we've checked in the proof of above Proposition), such operators preserve the space $\Psi^m(T_{t'})\otimes_{\Psi^m(\gL)}\gL$.  

(iii) If $\f$ is regular at $q=1/\z$, where $\z\neq 0,\infty$ is not a root of unity, $\g$ is obviously regular there too. \qed

\medskip

{\tt Corollary} (the $q$-string equation). {\em The range $\L \subset \K$ of permutation -equivariant J-function is invariant under the multiplication operators:
  \[ \f \mapsto e^{\textstyle \sum_{k>0}\Psi^k(\epsilon)/k(1-q^k)}\, \f,\ \ \ \epsilon \in \gL_{+}.\]}

{\tt Proof:} Use Theorem 1 with $D=1$.

\section*{Examples}

{\tt Example 1: $d=0$.}
In degree $0$, i.e. modulo Novikov's variables, the cone $\L \subset \K$ coincides with the cone $\L_{pt}$ over the $\lambda$-algebra $K^0(X)\otimes \gL$. Theorem 1 and Proposition allow one to recover the part of $\L_{pt}$ over the $\lambda$-algebra $\gL'=K^0(X)_{pr}\otimes \gL$, where by $K^0(X)_{pr}$ (the {\em primitive} part) we denote the part of the ring $K^0(X)$ generated by line bundles.

Let monomials $P^a:=P_1^{\a_1}\cdots P_K^{a_K}$ run a basis of
$K^0(X)_{pr}$. Applying the above theorem to the finite difference operator
\[ D=\sum_a \epsilon_a P^a q^{a Q\p_{Q}}:=\sum_a \epsilon_a \prod_{i=1}^K P_i^{a_i}q^{a_i Q_i\p_{Q_i}}, \ \ \epsilon_a \in \gL_{+},\] 
and acting on the point $\J \equiv 1-q $ {\em modulo Novikov's variables},
we recover over $\gL'$ the small J-function of the point: 
\[  (1-q) e^{\textstyle \sum_a \sum_{k>0} \Psi^k(\epsilon_a) P^{ka}/k(1-q^k)} \equiv
1-q+\sum_a \epsilon_a P^a \mod \K_{-}.\]
Furthermore, applying linear combinations
\[ \sum_a c_a(q) P^a q^{a Q\p_Q}\]
with coefficients $c_a \in \gL [q,q^{-1}]$ which are arbitrary Laurent polynomials in $q$,
we get, according to Proposition, points in the same ruling space of the cone $\L$. Modulo Novikov's variables this effectively results in multiplying by
arbitrary elements $\sum_a c_a(q) P^a$ from $\gL' [q,q^{-1}]$, and therefore yields the entire cone $\L_{pt}$ over $\gL'$.

\medskip

{\tt Example 2: $X=\CC P^1$.} We know\footnote{From various sources: Part IV (by localization), or \cite{GiTo} (by adelic characterization), or \cite{GiL} (by toric compactifications).} one point on $\L=\L_{\CC P^1}$, the small J-function:
\[ \J(0)=(1-q)\sum_{d\geq 0}\frac{Q^d}{(1-Pq)^2(1-Pq^2)^2\cdots (1-Pq^d)^2}.\]
Here $P=\O(-1)$ is the generator of $K^0(\CC P^1)$. It satisfies the relation $(1-P)^2=0$.  The K-theoretic Poincar\'e pairing is determined by
\[ \chi (\CC P^1; \phi(P)) = \Res_{P=1} \frac{\phi(P)}{(1-P)^2}\frac{dP}{P}.\]
We use Theorem 1 with the operator $D=\la + \ge P q^{Q\p_Q}$, $\la,\ge \in \gL_{+}$, and obtain a 2-parametric family of points on $\L_{\CC P^1}$: 
\begin{align*} &e^{\textstyle \sum_{k>0} (\Psi^k(\la)+\Psi^k(\ge) P^k q^{kQ\p_Q})/k(1-q^k)}\, \J(0) =\\  &(1-q) e^{\textstyle \sum_{k>0}\Psi^k(\la)/k(1-q^k)} \sum_{d\geq 0} \frac{Q^d\, e^{\textstyle \sum_{k>0}\Psi^k(\ge) P^k q^{kd}/k(1-q^k)}}{(1-Pq)^2(1-Pq^2)^2\cdots (1-Pq^d)^2}.\end{align*}

Examine now two specializations.

Firstly, as a consistency check, let us extract from this the {\em degree-$1$ part of $\F_0$.} Modulo $Q^2$, we are left with 
\[ \J \equiv e^{\sum_{k>0}\Psi^k(\la)/k(1-q^k)}\left(1-q+\frac{(1-q)Q}{(1-Pq)^2} \, e^{\sum_{k>0}\Psi^k(\ge)P^kq^k/k(1-q^k)}\right).\]
Modulo $\K_{-}$ (and $Q^2$), we have: $[\J]_{+}\equiv 1-q+\la+\ge P$.  
According to Part VII, Corollary 3,
\[ \F_0(\t) = -\frac{1}{2}\Omega ([\J(\t)]_{+},\J(\t))-\frac{1}{2}(\Psi^2(\t(1)),1).\]
For degree $d=1$ part $\J_1$ of $\J$, we have
\begin{align*} -\Omega ([\J]_{+}, \J_1) &= \Res_{q=0,\infty} \left(1-\frac{1}{q}+\la+\ge P, \frac{(1-q)}{(1-Pq)^2}\, e^{A(q)}\right) \frac{dq}{q}, \end{align*}
where $A(q)=\sum_{k>0} (\Psi^k(\la)+\Psi^k(\ge)P^kq^k)/k(1-q^k)$. The 1-form has no pole at $q=\infty$. Since $(\left(1-q)/(1-Pq)^2\right)'_{q=0}= 2P-1$, and $A'(0)=\la+\ge P$, the residue at $q=0$ is calculated as
\begin{align*} \left(1+\la+\ge P , e^{A(0)}\right) - \left(1,  (2P-1) e^{A(0)} + (\la+\ge P) e^{A(0)}\right) =\\
  \Res_{P=1} \frac{2(1-P) e^{A(0)}}{(1-P)^2}\frac{dP}{P} = 2 e^{A(0)}=2 e^{
    \sum_{k>0}\Psi^k(\la)/k} .\end{align*}
Let us check this rather trivial result ``by hands''. The degree $d=1$ part of $\F_0(t)$ at $t=\la+\ge P$ is defined as
$\sum_{n\geq 0} \lan \la+\ge P, \dots, \la+\ge P\ran_{0,n,1}^{S_n}$. Since there is only one rational curve of degree $1$ in $\CC P^1$, the moduli space
$X_{0,n,1}=\M_{0,n}(\CC P^1, 1)$ is obtained from $(\CC P^1)^n$ by some blow-ups along the diagonals. The evaluation maps $\ev_i: X_{0,n,1}\to \CC P^1$ factor through $(\CC P^1)^n$ as the projections  $(\CC P^1)^n\to \CC P^1$. Therefore the correlator sum can be evaluated as
\begin{align*} \sum_{n\geq 0} \left(H^*\left(\CC P^1; \la +\ge P\right)^{\otimes n}\right)^{S_n} = \sum_{n\geq 0} (\la^{\otimes n})^{S_n} \end{align*}
because for $P=\O(-1)$ we have $H^*(\CC P^1; P)=0$. Let us remind from Part I that for elements of a $\la$-algebra, 
\[ (\la^{\otimes n})^{S_n} :=\frac{1}{n!}\sum_{h\in S_n} \prod_{k>0}  \Psi^{l_k(h)}(\la),\]
where $l_k(h)$ is the number of cycles of length $k$ in the permutation $h$.
Thus, the correlator sum indeed coincides with $e^{\sum_{k>0} \Psi^k(\la)/k}$. 

\pagebreak

Secondly, let us return to our 2-parametric family of points on $\L_{\CC P^1}$, and specialize it to the {\em symmetrized} theory, where only the $S_n$-invariant part of sheaf cohomology is taken into account. For this, we specialize the $\la$-algebra to $\gL = \QQ[[\la, \ge, Q]]$ with  $\Psi^k(\la)=\la^k, \Psi^k(\ge)=\ge^k$ (and $\Psi^k(Q)=Q^k$ as before).  Some simplifications ensue. Since
\[ q^{kd}=1-(1-q^k)(1+q^k+\cdots+q^{k(d-1)}),\]
we have
\begin{align*}  e^{\textstyle \sum_{k>0}\ge^kP^kq^{kd}/k(1-q^k)} &= e^{\textstyle \sum_{k>0}\ge^kP^k/k(1-q^k)} \prod_{r=0}^{d-1} e^{-\sum_{k>0} \ge^kP^kq^{rk}/k} \\ &=e^{\textstyle \sum_{k>0}\ge^kP^k/k(1-q^k)} \prod_{r=0}^{d-1}(1-\ge Pq^r).\end{align*}    
Thus, we obtain the following 2-parametric family of points on $\L^{sym}_{\CC P^1}$:
\[ \J^{sym}_{\CC P^1}=(1-q)e^{\textstyle \sum_{k>0}(\la^k+\ge^kP^k)/k(1-q^k)}
 \sum_{d\geq 0}Q^d\,\frac{\prod_{r=0}^{d-1}(1-\ge P q^r)}{\prod_{r=1}^d(1-Pq^r)^2}.\]
Note that the projection of this series to $\K_{+}$ along $\K_{-}$ picks contributions only from the terms with $d=0$ and $k=1$: 
\[ [\J^{sym}_{\CC P^1}]_{+} = 1-q+\la +\ge P.\]
Therefore the series represents {\em the small J-function of the symmetrized quantum K-theory of $\CC P^1$}. The exponential factor is actually equal to $\exp_q(\la/(1-q))\exp_q(\ge P/(1-q))$. Thus, we obtain:
\begin{align*} &\J^{sym}_{\CC P^1}(\la+\ge P) =_{\mod (1-P)^2}  \\
 &(1-q)\sum_{m=0}^{\infty}\sum_{l=0}^{\infty}\sum_{d=0}^{\infty}\frac{\la^m\, \ge^l\, P^l\, Q^d \, \prod_{r=0}^{d-1}(1-\ge Pq^r)}{\prod_{t=1}^m(1-q^t) \prod_{s=1}^l(1-q^s)\prod_{r=1}^d(1-Pq^r)} . \end{align*}

\medskip

\section*{Reconstruction theorems}

As in Example 1, assume that $p_1,\dots, p_K$ is a numerically effective integer basis in $H^2(X, \QQ)$, that Novikov's monomials $Q^d=Q_1^{d_1}\dots Q_K^{d_K}$ represent degree $d$ holomorphic curves in $X$ in coordinates $d_i=p_i(d)$ on $H_2(X)$, that $P_i$ are line bundles with $\ch P_i=e^{-p_i}$, and that monomials $P^a=P_1^{a_1}\cdots P_K^{a_k}$ run a basis in $K^0(X)_{pr}$, the primitive part of the K-ring. We also write $a.d$ for the value  $\sum a_id_i$ of $-c_1(P^a)$ on $d$.  

\pagebreak

{\tt Theorem 2} (explicit reconstruction). {\em Let $I=\sum_d I_d Q^d$ be a point in the range $\L \subset \K$ of the J-function of permutation-equivariant quantum K-theory on $X$, written as a vector-valued series in Novikov's variables. Then the following family also lies in $\L$:  
  \[ \sum_d I_d Q^d e^{\textstyle \sum_{k>0} \sum_a\Psi^k(\ge_a) P^{ka} q^{k(a.d)}/k(1-q^k)}, \ \ \ge_a \in \gL_{+} .\]
Moreover, for arbitrary Laurent polynomials $c_a\in \gL [q,q^{-1}]$, the following series also lies in $\L$:
\[  \sum_d I_d Q^d e^{\textstyle \sum_{k>0} \sum_a\Psi^k(\ge_a) P^{ka} q^{k(a.d)}/k(1-q^k)}\sum_a c_a(q) P^aq^{a.d}. \]
Furthermore, when $K^0(X)=K^0(X)_{pr}$, the whole cone $\L\subset \K$ is parameterized this way.}  

\medskip

{\tt Proof.} We first work over $\gL'$ freely generated as $\lambda$-algebra by
the ``time'' variables $\ge_a$, and use Theorem 1 with the $Q$-independent
finite difference operator $D = \sum_a \ge_a P^a q^{aQ\p_Q}$. We conclude that the family
\[ e^{\sum_{k>0} \sum_a\Psi^k(\ge_a) P^{ka} q^{ka Q\p_Q}/k(1-q^k)} I=
  \sum_d I_d Q^d e^{\sum_{k>0}\sum_a \Psi^k(\ge_a)P^aq^{k(a.d)}/k(1-q^k)} \]
lies in the cone $\L$, defined over $\gL'$. 
To obtain the second statement, we apply Proposition, using finite difference operators $\sum_a c_a(q) P_a q^{a Q\p_Q}$.
Afterwards we specialize the ``times'' $\ge_a$ to any values $\ge_a \in \gL$ (which at this point may become dependent on $Q$). Finally, when $K^0(X)=K^0(X)_{pr}$, we use the formal Implicit Function Theorem to conclude that the whole cone $\L$ is parameterized, because this is true modulo Novikov's variables, as Example 1 shows. \qed
  
\medskip

{\tt Example: $X=\CC P^N$.} According to Theorem 2, the entire cone $\L$ is parameterized as follows:
\[ \J = (1-q)\sum_{d\geq 0}Q^d\,\frac{e^{\sum_{k>0}\sum_{a=0}^N \Psi^k(\ge_a) P^{ka}q^{kad}/k(1-q^k)} \sum_{a=0}^N c_a(q) P^aq^{ad}}{(1-Pq)^{N+1}(1-Pq^2)^{N+1}\cdots (1-Pq^d)^{N+1}}.\]
Of course, this is obtained by applying Theorem 2 to
the small J-function $\J(0)$ from \cite{GiL} (also \cite{GiTo},  or
Parts II--IV in the non-equivariant limit).
Here $\ge_a\in \gL$, $c_a(q)$ are arbitrary Laurent polynomials in $q$ with coefficients in $\gL$, and $P^a$, $a=0,\dots, N$, $P=\O(-1)$, are used for a basis in $K^0(X)$. Perhaps, the basis $(1-P)^a, a=0,\dots, N$, is more useful (cf. \cite{GiE}), and we get yet another parameterization of $\L$:
\[ (1-q)\sum_{d\geq 0}Q^d\,\frac{e^{\sum_{k>0}\sum_{a=0}^N \Psi^k(\ge_a) (1-P^kq^{kd})^a/k(1-q^k)} \sum_{a=0}^N c_a(q) (1-Pq^{d})^a} {(1-Pq)^{N+1}(1-Pq^2)^{N+1}\cdots (1-Pq^d)^{N+1}}.\]

\medskip

We return now to the context of Part VII, where we studied the mixed J-function $\J(\x,\t)$, involving two types of inputs: permutable $\t$ and non-permutable $\x$, both taken from $\K_{+}$. The cone $\L\subset \K$ represents the range of $\t\mapsto \J(0,\t)$. Recall that according to the general theory, it is the union of {\em ruling spaces} $(1-q) T_t$, where $t=\T(\t)$ is given by a certain non-linear map
\[ \T: K^0(X)\otimes \gL_{+}\oplus (1-q)\K_{+}\to K^0(X)\otimes \gL_{+}.\]
At the same time, for a fixed value of $\t$, the range of the {\em ordinary} J-function $\x\mapsto \J(\x,\t)$ is an overruled Lagrangian cone $\L_t\subset \K$, which shares with $\L$ one ruling space, $T_t$, corresponding to $t=\T(\t)$.
Each tangent space of each cone  $\L_t$ is tangent to $\L_t$ along one of the ruling spaces (e.g. $T_t$ is tangent along $(1-q)T_t$), and is related with this ruling space by the multiplication by $1-q$. As a consequence, not only each ruling (and tangent) space of each $\L_t$ is a $D_q$-module (which is proved on the basis of adelic characterization as in Proposition above), but also each cone $\L_t$ is invariant under the flow
\[ \f \mapsto e^{\ge D(Q, Pq^{Q\p_Q}, q)/(1-q)}\f ,\]
where $D \in \D_q$. We use this to reconstruct the family $\L_t$. 

\medskip

{\tt Theorem 3.} {\em Let $I=\sum I_d Q^d$ (as in Theorem 2). Then   
  \[ I(\ge)=\sum_d I_d(\ge) Q^d:=\sum_d I_d Q^d e^{\sum_{k>0} \sum_a\Psi^k(\ge_a) P^{ka} q^{k(a.d)}/k(1-q^k)}, \ \ \ge_a \in \gL_{+} \]
  represent a family of points on the cones $\L_{t(\ge)}$ (one point on each cone), and the following family of points, parameterized by $\tau_a\in \gL$ and by $c_a\in \gL[q,q^{-1}]$, lies on $\L_{t(\ge)}$:
  \[ \sum_d I_d(\ge)\, Q^d \, e^{\textstyle \sum_a \tau_a P^aq^{a.d}/(1-q)} 
  \sum_a c_a(q) P^aq^{a.d}. \]
  Moreover, if $K^0(X)=K^0(X)_{pr}$, for each $t\in K^0(X)\otimes \gL_{+}$ the whole cone $\L_t$ is thus parameterized.}

\medskip

{\tt Proof.} It is clear from computation modulo Novikov's variables that the family $I(\ge)$ has no tangency with the ruling spaces, hence represents at most one point from each $\L_t$ (and does represent one, when $K^0(X)=K^0(X)_{pr}$). Given one point, $I(\ge)$, on $\L_{t(\ge)}$, we generate more points by machinery discussed above: applying the commuting flows
\[ e^{\sum_a \tau_a P^aq^{O\p_Q}/(1-q)} I(\ge)=\sum_d Q^d I_d(\ge)\, e^{\sum_a \tau_a P^a q^{a.d}/(1-q)},\]
followed by the application of the operators $\sum_a c_a(q) P^q q^{aQ\p_Q}$, where $\tau_a$ and the coefficients of $c_a$ are independent variables.  Afterwards they can be specialized to some values in $\gL$ (in particular, depending on $Q$). In the case when $K^0(X)=K^0(X)_{pr}$, it follows from the Implicit Function Theorem and Example 1 about the limit to $d=0$, that the entire cone $\L_t$ for each $t$ is thus obtained. \qed

\medskip

{\tt Example: $X=\CC P^N$.} It follows that for fixed values of $\ge_a$, the
corresponding cone $\L_{t(\ge)}$ is parameterized as
\[ (1-q)\sum_{d\geq 0}Q^d\,\frac{e^{\sum_{a=0}^N \left(\tau_a\frac{P^aq^{ad}}{1-q}+\sum_{k>0}\Psi^k(\ge_a) \frac{P^{ka}q^{kad}}{k(1-q^k)}\right)}\sum_{a=0}^N c_a(q) P^aq^{ad}}{(1-Pq)^{N+1}(1-Pq^2)^{N+1}\cdots (1-Pq^d)^{N+1}},\]
and all $\L_t$ are so obtained.

\medskip

{\tt Remarks.} Reconstruction theorems in quantum cohomology and K-theory go back to Kontsevich--Manin \cite{KM} and Lee--Pandharipande \cite{LP} respectively. Theorem 3 is a slight generalization (from the case $t=0$) of the ``explicit reconstruction'' result \cite{GiE} in the ordinary (non-permutation-equivariant) quantum K-theory, which in its turn mimics the results of quantum cohomology theory already found in \cite{CF-K, I}, and shares the methods based on finite difference operators with the K-theoretic results of \cite{IMT}.

Theorems of this section show that when $K^0(X)$ is generated by line bundles, the entire range $\L$ of the J-function in the permutation-equivariant genus-0 quantum K-theory of $X$, as well as the entire family $\L_t$ of the overruled Lagrangian cones representing the 
``ordinary'' J-functions, depending on the permutable parameter, $t$, can be explicitly represented in a parametric form, {\em given one point on any of these cones}. In essence, all genus-0 K-theoretic GW-invariants of $X$, permutation-equivariant, ordinary, or mixed, are thereby reconstructed from any one point: a $K^0(X)$-valued series $\sum_d I_d Q^d$ in Novikov's variables. 

In the case of a toric $X$, the results of Part V exhibit such a point in the form of the {\em $q$-hypergeometric series} mirror-symmetric to $X$. Needless to say, the same applies to toric bundles spaces, or super-bundles (a.k.a. toric complete intersections), as well as to the torus-equivariant versions of K-theoretic GW-invariants. Thus ``all'' (torus-equivariant or not; permutation-equivariant, ordinary, or mixed) K-theoretic genus-0 GW-invariants of toric manifolds, toric bundles, or toric complete intersections are computed in a geometrically explicit form, illustrated by the above example.    


\enddocument